# ABOUT THE CHARACTERISTIC FUNCTION OF A SET

Prof. Mihály Bencze, Department of Mathematics,
University of Braşov, Romania

Prof. Florentin Smarandache, Chair of Department of Math & Sciences, University of New Mexico, 200 College Road, Gallup, NM 87301, USA, E-mail: smarand@unm.edu

**Abstract**:
In this paper we give a method, based on the characteristic function of a set, to solve some difficult problems of set theory found in undergraduate studies.

**Definition**: Let's consider $A \subset E \neq \varnothing$ (a universal set), then $f_A : E \to \{0, 1\}$, where the function $f_A(x) = \begin{cases} 1, & \text{if } x \in A \\ 0, & \text{if } x \notin A \end{cases}$ is called the characteristic function of the set $A$.

**Theorem 1:** Let's consider $A, B \subset E$. In this case $f_A = f_B$ if and only if $A = B$.

**Proof.**

$$f_A(x) = \begin{cases} 1, & \text{if } x \in A = B \\ 0, & \text{if } x \notin A = B \end{cases} = f_B(x)$$

Reciprocally: For any $x \in A$, $f_A(x) = 1$, but $f_A = f_B$, therefore $f_B(x) = 1$, namely $x \in B$ from where $A \subset B$. The same way we prove that $B \subset A$, namely $A = B$.

**Theorem 2:** $f_{\tilde{A}} = 1 - f_A$, $\tilde{A} = C_E A$.

**Prof.**

$$f_{\tilde{A}}(x) = \begin{cases} 1, & \text{if } x \in \tilde{A} \\ 0, & \text{if } x \notin \tilde{A} \end{cases} = \begin{cases} 1, & \text{if } x \notin A \\ 0, & \text{if } x \in A \end{cases} = \begin{cases} 1-0, & \text{if } x \notin A \\ 1-1, & \text{if } x \in A \end{cases} = 1 - \begin{cases} 0, & \text{if } x \notin A \\ 1, & \text{if } x \in A \end{cases} = 1 - f_A(x)$$

**Theorem 3:** $f_{A \cap B} = f_A * f_B$.

**Proof.**

$$f_{A\cap B}(x) = \begin{cases} 1, & \text{if } x \in A \cap B \\ 0, & \text{if } x \notin A \cap B \end{cases} = \begin{cases} 1, & \text{if } x \in A \text{ and } x \in B \\ 0, & \text{if } x \notin A \text{ or } x \notin B \end{cases} = \begin{cases} 1, & \text{if } x \in A, x \in B \\ 0, & \text{if } x \in A, x \notin B \\ 0, & \text{if } x \notin A, x \in B \\ 0, & \text{if } x \notin A, x \notin B \end{cases} =$$

$$= \left( \begin{cases} 1 & \text{if } x \in A \\ 0 & \text{if } x \notin A \end{cases} \right) \left( \begin{cases} 1 & \text{if } x \in B \\ 0 & \text{if } x \notin B \end{cases} \right) = f_A(x) f_B(x).$$

The theorem can be generalized by induction:

**Theorem 4:** $f_{\bigcap_{k=1}^{n} A_k} = \prod_{k=1}^{n} f_{A_k}$

*Consequence.* For any $n \in \mathbb{N}^*$, $f_M^n = f_M$.

**Proof.** In the previous theorem we chose $A_1 = A_2 = \ldots = A_n = M$.

**Theorem 5:** $f_{A \cup B} = f_A + f_B - f_A f_B$.

**Proof.**
$$f_{A \cup B} = f_{\overline{\overline{A \cup B}}} = f_{\overline{\overline{A} \cap \overline{B}}} = 1 - f_{\overline{A} \cap \overline{B}} = 1 - f_{\overline{A}} f_{\overline{B}} = 1 - (1 - f_A)(1 - f_B) = f_A + f_B - f_A f_B$$

It can be generalized by induction:

**Theorem 6:** $f_{\bigcup_{k=1}^{n} A_k} = \sum_{k=1}^{n} (-1)^{k-1} \sum_{1 \le i_1 < \ldots < i_k \le n} (-1)^{k-1} f_{A_{i_1}} f_{A_{i_2}} \cdots f_{A_{i_k}}$

**Theorem 7:** $f_{A-B} = f_A(1 - f_B)$

**Proof.** $f_{A-B} = f_{A \cap \overline{B}} = f_A f_{\overline{B}} = f_A(1 - f_B)$.

It can be generalized by induction:

**Theorem 8:** $f_{A_1 - A_2 - \ldots - A_n} = \sum_{k=1}^{n} (-1)^{k-1} f_{A_{i_1}} f_{A_{i_2}} \cdots f_{A_{i_k}}$.

**Theorem 9:** $f_{A \Delta B} = f_A + f_B - 2 f_A f_B$

**Proof.**
$$f_{A \Delta B} = f_{A \cup B - A \cap B} = f_{A \cup B}(1 - f_{A \cap B}) = (f_A + f_B - f_A f_B)(1 - f_A f_B) = f_A + f_B - 2 f_A f_B.$$

It can be generalized by induction:

**Theorem 10:** $F_{\Delta_{k=1}^{n} A_k} = \sum_{k=1}^{n} (-2)^{k-1} \sum_{1 \le i_1 < \ldots < i_k \le n} f_{A_{i_1} A_{i_2} \ldots A_{i_k}}$.

**Theorem 11:** $f_{A \times B}(x, y) = f_A(x) f_B(y)$.

**Proof.** If $(x, y) \in A \times B$, then $f_{A \times B}(x, y) = 1$ and $x \in A$, namely $f_A(x) = 1$ and $y \in B$, namely $f_B(y) = 1$, therefore $f_A(x) f_B(y) = 1$. If $(x, y) \notin A \times B$, then $f_{A \times B}(x, y) = 0$ and $x \notin A$, namely $f_A(x) = 0$ or $y \notin B$, namely $f_B(y) = 0$, therefore $f_A(x) f_B(y) = 0$.

This theorem can be generalized by induction.

**Theorem 12:** $f_{\times_{k=1}^{n} A_k}(x_1, x_2, \ldots, x_n) = \prod_{k=1}^{n} f_{A_k}(x_k)$.

**Theorem 13**: (De Morgan) $\overline{\bigcup_{k=1}^{n} A_k} = \bigcap_{k=1}^{n} \overline{A_k}$

**Proof.**

$$f_{\overline{\bigcup_{k=1}^{n} A_k}} = 1 - f_{\bigcup_{k=1}^{n} A_k} = 1 - \sum_{k=1}^{n}(-1)^{k-1} \sum_{1 \leq i_1 < \ldots < i_k \leq n} f_{A_{i_1}} f_{A_{i_2}} \cdots f_{A_{i_k}} = \prod_{k=1}^{n}(1 - f_{A_k}) = \prod_{k=1}^{n} f_{\overline{A_k}} = f_{\bigcap_{k=1}^{n} \overline{A_k}}.$$

We prove in the same way the following theorem:

**Theorem 14:** (De Morgan) $\overline{\bigcap_{k=1}^{n} A_k} = \bigcup_{k=1}^{n} \overline{A_k}$.

**Theorem 15:** $\left(\bigcup_{k=1}^{n} A_k\right) \cap M = \bigcup_{k=1}^{n}(A_k \cap M)$.

**Proof.**

$$f_{\left(\bigcup_{k=1}^{n} A_k\right) \cap M} = f_{\bigcup_{k=1}^{n} A_k} f_M = \sum_{k=1}^{n}(-1)^{k-1} \sum_{1 \leq i_1 < \ldots < i_k \leq n} f_{A_{i_1}} f_{A_{i_2}} \cdots f_{A_{i_k}} f_M = \sum_{k=1}^{n}(-1)^{k-1} \sum_{1 \leq i_1 < \ldots < i_k \leq n} f_{A_{i_1}} f_{A_{i_2}} \cdots f_{A_{i_k}} f_M^k =$$

$$= \sum_{k=1}^{n}(-1)^{k-1} \sum_{1 \leq i_1 < \ldots < i_k \leq n} f_{A_{i_1} \cap M} f_{A_{i_2} \cap M} \cdots f_{A_{i_k} \cap M} = f_{\bigcup_{k=1}^{n}(A_k \cap M)}$$

In the same way we prove that:

**Theorem 16:** $\left(\bigcap_{k=1}^{n} A_k\right) \cup M = \bigcap_{k=1}^{n}(A_k \cup M)$.

**Theorem 17:** $\left(\Delta_{k=1}^{n} A_k\right) \cap M = \Delta_{k=1}^{n}(A_k \cap M)$

**Application.**
$\left(\Delta_{k=1}^{n} A_k\right) \cup M = \Delta_{k=1}^{n}(A_k \cup M)$ if and only if $M = \Phi$.

**Theorem 18:** $M \times \left(\bigcup_{k=1}^{n} A_k\right) = \bigcup_{k=1}^{n}(M \times A_k)$

**Proof.**

$$f_{M\times\left(\bigcup_{k=1}^{n}A_k\right)}(x,y)=f_M(y)f_{\bigcup_{k=1}^{n}A_k}(x)=\sum_{k=1}^{n}(-1)^{k-1}\sum_{1\leq i_1<...<i_k\leq n}f_{A_{i_1}}(x)f_{A_{i_2}}(x)...f_{A_{i_k}}(x)f_M(y)=$$

$$=\sum_{k=1}^{n}(-1)^{k-1}\sum_{1\leq i_1<...<i_k\leq n}f_{A_{i_1}}(x)f_{A_{i_2}}(x)...f_{A_{i_k}}(x)f_M^k(y)=$$

$$=\sum_{k=1}^{n}(-1)^{k-1}\sum_{1\leq i_1<...<i_k\leq n}f_{A_{i_1}\times M}(x,y)...f_{A_{i_k}\times M}(x,y)=f_{\bigcup_{k=1}^{n}(M\times A_k)}$$

In the same way we prove that:

**Theorem 19:** $M\times\left(\bigcap_{k=1}^{n}A_k\right)=\bigcap_{k=1}^{n}(M\times A_k)$.

**Theorem 20:** $M\times(A_1-A_2-...-A_n)=(M\times A_1)-(M\times A_2)-...-(M\times A_n)$.

**Theorem 21:** $(A_1-A_2)\cup(A_2-A_3)\cup...\cup(A_{n-1}-A_n)\cup(A_n-A_1)=\bigcup_{k=1}^{n}A_k-\bigcap_{k=1}^{n}A_k$

**Proof 1.**

$$f_{(A_1-A_2)\cup...\cup(A_n-A_1)}=\sum_{k=1}^{n}(-1)^{k-1}\sum_{1\leq i_1<...<i_k\leq n}f_{A_{i_1}-A_{i_2}}...f_{A_{i_k}-A_{i_1}}=$$

$$=\sum_{k=1}^{n}(-1)^{k-1}\sum_{1\leq i_1<...<i_k\leq n}(f_{A_{i_1}}-f_{A_{i_1}}f_{A_{i_2}})...(f_{A_{i_k}}-f_{A_{i_k}}f_{A_{i_1}})=$$

$$=\sum_{k=1}^{n}(-1)^{k-1}\sum_{1\leq i_1<...<i_k\leq n}f_{A_{i_1}}...f_{A_{i_k}}\left(1-\prod_{p=1}^{n}f_{A_p}\right)=f_{\bigcup_{k=1}^{n}A_k}\left(1-f_{\bigcap_{k=1}^{n}A_k}\right)=f_{\bigcup_{k=1}^{n}A_k-\bigcap_{k=1}^{n}A_k}.$$

**Proof 2.** Let's consider $x\in\bigcup_{i=1}^{n}(A_i-A_{i+1})$, (where $A_{n+1}=A_1$), then there exists $k$ such that $x\in(A_k-A_{k+1})$, namely $x\notin(A_k\cap A_{k+1})\subset A_1\cap A_2\cap...\cap A_n$, namely $x\notin A_1\cap A_2\cap...\cap A_n$, and $x\in\bigcup_{k=1}^{n}A_k-\bigcap_{k=1}^{n}A_k$.

Now we prove the inverse statement:

Let's consider $x\in\bigcup_{k=1}^{n}A_k-\bigcap_{k=1}^{n}A_k$, we show that there exists $k$ such that $x\in A_k$ and $x\notin A_{k+1}$. On the contrary, it would result that for any $k\in\{1,2,...,n\}$, $x\in A_k$ and $x\in A_{k+1}$ namely $x\in\bigcup_{k=1}^{n}A_k$, it results that there exists $p$ such that $x\in A_p$, but from the previous reasoning it results that $x\in A_{p+1}$, and using this we consequently obtain that $x\in A_k$ for $k=\overline{p,n}$. But from $x\in A_n$ we obtain that $x\in A_1$, therefore, it results that $x\in A_k$, $k=\overline{1,p}$, from where $x\in A_k$, $k=\overline{1,n}$, namely $x\in A_1\cap...\cap A_n$, that is a contradiction. Thus there exists $r$ such that $x\in A_r$ and $x\notin A_{r+1}$, namely $x\in(A_r-A_{r+1})$ and therefore $x\in\bigcup_{k=1}^{n}(A_k-A_{k+1})$.

In the same way we prove the following theorem:

**Theorem 22:** $(A_1 \Delta A_2) \cup (A_2 \Delta A_3) \cup ... \cup (A_{n-1} \Delta A_n) = \bigcup_{k=1}^{n} A_k - \bigcap_{k=1}^{n} A_k$.

**Theorem 23:**
$(A_1 \times A_2 \times ... \times A_k) \cap (A_{k+1} \times A_{k+2} \times ... \times A_{2k}) \cap (A_n \times A_1 \times ... \times A_{k-1}) = (A_1 \cap A_2 \cap ... \cap A_n)^k$.

**Proof.** $f_{(A_1 \times ... \times A_k) \cap ... \cap (A_n \times A_1 \times ... \times A_{k-1})}(x_1,...,x_n) =$

$= f_{A_1 \times ... \times A_k}(x_1,...,x_n) ... f_{A_n \times ... \times A_{k-1}}(x_1,...,x_n) =$

$= (f_{A_1}(x_1)...f_{A_k}(x_k))...(f_{A_n}(x_n)...f_{A_{k-1}}(x_{k-1})) =$

$= f_{A_1}^k(x_1)...f_{A_n}^k(x_n) = f_{A_1 \cap ... \cap A_n}^k(x_1,...,x_n) =$

$= f_{(A_1 \cap ... \cap A_n)^k}(x_1,...,x_n)$.

**Theorem 24.** $(P(E), \cup)$ is a commutative monoid.

**Proof.** For any $A, B \in P(E)$; $A \cup B \in P(E)$, namely the intern operation. Because $(A \cup B) \cup C = A \cup (B \cup C)$ is associative, $A \cup B = B \cup A$ commutative, and because $A \cup \emptyset = A$ then $\emptyset$ is the neutral element.

**Theorem 25:** $(P(E), \cap)$ is a commutative monoid.

**Proof.** For any $A, B \in P(E)$; $A \cap B \in P(E)$ namely intern operation. $(A \cap B) \cap C = A \cap (B \cap C)$ associative, $A \cap B = B \cap A$, commutative $A \cap E = A$, $E$ is the neutral element.

**Theorem 26:** $(P(E), \Delta)$ is an abelian group.

**Proof.** For any $A, B \in P(E)$; $A \Delta B \in P(E)$, namely the intern operation. $A \Delta B = B \Delta A$ commutative. The proof of associativity is in the XII$^{th}$ grade manual as a problem. We'll prove it using the characteristic function of the set.

$f_{(A \Delta B) \Delta C} = 4 f_A f_B f_C - 2 f_A f_B + f_B f_C + f_C f_A + f_A + f_B + f_C = f_{A \Delta (B \Delta C)}$.

Because $A \Delta \emptyset = A$, $\emptyset$ is the neutral element and because $A \Delta A = \emptyset$; the symmetric element of $A$ is $A$ itself.

**Theorem 27:** $(P(E), \Delta, \cap)$ is a commutative Boole ring with a divisor of zero.

**Proof.** Because the previous theorem satisfies the commutative ring axioms, the first part of the theorem is proved. Now we prove that it has a divisor of zero. If $A \neq \emptyset$ and $B \neq \emptyset$ are two disjoint sets, then $A \cap B = \emptyset$, thus it has divisor of zero. From Theorem 17 we get that it is distributive for $n = 2$. Because for any $A \in P(E)$; $A \cap A = A$ and $A \Delta A = \emptyset$ it also satisfies the Boole-type axioms.

**Theorem 28:** Let's consider $H = \{f \mid f : E \to \{0,1\}\}$, then $(H, \oplus)$ is an abelian group, where $f_A \oplus f_B = f_A + f_B - 2f_A f_B$ and $(P(E), \Delta) \cong (H, \oplus)$.

**Proof.** Let's consider $F : P(E) \to H$, where $f(A) = f_A$, then, from the previous theorem we get that it is bijective and because $F(A\Delta B) = f_{A\Delta B} = F(A) \oplus F(B)$ it is compatible.

**Theorem 29:** $card(A_1 \Delta A_n) \leq card(A_1 \Delta A_2) + card(A_2 \Delta A_3) + \ldots + card(A_{n-1} \Delta A_n)$.

**Proof.** By induction. If $n = 2$, then it is true, we show that for $n = 3$ it is also true. Because $(A_1 \cap A_2) \cup (A_2 \cap A_3) \subseteq A_2 \cup (A_1 \cap A_3)$;

$card((A_1 \cap A_2) \cup (A_2 \cap A_3)) \leq card(A_2 \cup (A_1 \cap A_3))$ but

$card(M \cup N) = cardM + cardN - card(M \cap N)$, and thus

$cardA_2 + card(A_1 \cap A_3) - card(A_1 \cap A_2) - card(A_2 \cap A_3) \geq 0$, can be written as

$cardA_1 + cardA_3 - 2card(A_1 \cap A_3) \leq$

$\leq (cardA_1 + cardA_2 - 2card(A_1 \cap A_2)) + (cardA_2 + cardA_3 - 2card(A_2 \cap A_3))$.

But because of

$(M \Delta N) = cardM + cardN - 2card(M \cap N)$

then $card(A_1 \Delta A_3) \leq card(A_1 \Delta A_2) + card(A_2 \Delta A_3)$. The proof of this step of the induction relies on the above method.

**Theorem 30:** $(P^2(E), card(A\Delta B))$ is a metric space.

**Proof.** Let $d(A, B) = card(A \Delta B) : P(E) \times P(E) \to \Box$

1. $d(A, B) = 0 \Leftrightarrow card(A\Delta B) = 0 \Leftrightarrow card((A - B) \cup (B - A)) = 0$ but because $(A - B) \cap (B - A) = \emptyset$ we obtain $(A - B) + card(B - A) = 0$ and because $(A - B) = 0$ and $card(B - A) = 0$, then $A - B = \emptyset$, $B - A = \emptyset$, and $A = B$.

2. $d(A, B) = d(B, A)$ results from $A\Delta B = B\Delta A$.

3. As a consequence of the previous theorem $d(A, C) \leq d(A, B) + d(B, C)$.

As a result of the above three properties it is a metric space.

## PROBLEMS

**Problem 1.**

Let's consider $A = B \cup C$ and $f : P(A) \to P(A) \times P(A)$, where $f(x) = (X \cup B, X \cup C)$. Prove that $f$ is injective if and only if $B \cap C = \emptyset$.

*Solution 1.* If $f$ is injective. Then

$f(\emptyset) = (\emptyset \cup B, \emptyset \cup C) = (B, C) = ((B \cap C) \cup B, (B \cap C) \cup C) = f(B \cap C)$ from which we obtain $B \cap C = \emptyset$. Now reciprocally: Let's consider $B \cap C = \emptyset$, then $f(X) = f(Y)$; it results that $X \cup B = Y \cup B$ and $X \cup C = Y \cup C$ or $X = X \cup \varnothing = X \cup (B \cap C) = (X \cup B) \cap (X \cup C) = (Y \cup B) \cap (Y \cup C) = Y \cup (B \cap C) = Y \cup \varnothing = Y$
namely it is injective.

*Solution 2.* Let's consider $B \cap C = \varnothing$ passing over the set function $f(X) = f(Y)$ if and only if $X \cup B = Y \cup B$ and $X \cup C = Y \cup C$, namely $f_{X \cup B} = f_{Y \cup B}$ and $f_{X \cup C} = f_{Y \cup C}$ or $f_X + f_B - f_X f_B = f_Y + f_B - f_Y f_B$ and $f_X + f_C - f_X f_C = f_Y + f_C - f_Y f_C$ from which we obtain $(f_X - f_Y)(f_B - f_C) = 0$.

Because $A = B \cup C$ and $B \cap C = \varnothing$, we have

$$(f_B - f_C)(u) = \begin{cases} 1, & \text{if } u \in B \\ -1, & \text{if } u \in C \end{cases} \neq 0$$

therefore $f_X - f_Y = 0$, namely $X = Y$ and thus it is injective.

**Generalization.** Let $M = \bigcup_{k=1}^{n} A_k$ and $f : P(A) \to P^n(A)$, where

$f(X) = (X \cup A_1, X \cup A_2, ..., X \cup A_n)$.

Prove that $f$ is injective if and only if $A_1 \cap A_2 \cap ... \cap A_n = \varnothing$.

**Problem 2.** Let $E \neq \varnothing$, $A \in P(E)$, and $f : P(E) \to P(E) \times P(E)$, where $f(X) = (X \cap A, X \cup A)$.
  a. Prove that $f$ is injective
  b. Prove that $\{f(x), x \in P(E)\} = \{(M, N) \mid M \subset A \subset N \subset E\} = K$.
  c. Let $g : P(E) \to K$, where $g(X) = f(X)$. Prove that $g$ is bijective and compute its inverse.

*Solution.*
  a. $f(X) = f(Y)$, namely $(X \cap A, X \cup A) = (Y \cap A, Y \cup A)$ and then $X \cap A = Y \cap A$, $X \cup A = Y \cup A$, from where $X \Delta A = Y \Delta A$ or $(X \Delta A) \Delta A = (Y \Delta A) \Delta A$, $X \Delta (A \Delta A) = Y \Delta (A \Delta A)$, $X \Delta \varnothing = Y \Delta \varnothing$ and thus $X = Y$, namely $f$ is injective.

  b. $\{f(X), X \in P(E)\} = f(P(E))$. We'll show that $f(P(E)) \subset K$. For any $(M, N) \in f(P(E))$, $\exists X \in P(E) : f(X) = (M, N)$; $(X \cap A, X \cup A) = (M, N)$.
  From here $X \cap A = M$, $X \cup A = N$, namely $M \subset A$ and $A \subset N$
thus $M \subset A \subset N$, and, therefore $(M, N) \in X$.

  Now, we'll show that $K \subset f(P(E))$, for any $(M, N) \in K$, $\exists X \in P(E)$ such that $f(X) = (M, N)$. $f(X) = (M, N)$, namely $(X \cap A, X \cup A) = (M, N)$ from where $X \cap A = M$ and $X \cup A = N$, namely $X \Delta A = N - M$, $(X \Delta A) \Delta A = (N - M) \Delta A$, $X \Delta \varnothing = (N - M) \Delta A$,
$X = (N - M) \Delta A$, $X = (N \cap \overline{M}) \Delta A$,
$X = \big((N \cap \overline{M}) - A\big) \cup \big(A - (N \cap \overline{M})\big) = \big((N \cap \overline{M}) \cap \overline{A}\big) \cup \big(A \cap \overline{(N \cap \overline{M})}\big) =$
$= \big(N \cap (\overline{M} \cap \overline{A})\big) \cup \big(A \cap (\overline{N} \cap \overline{M})\big) = (N \cap \overline{A}) \cup \big((A \cap \overline{N}) \cup (A \cap M)\big) =$

$= (N \cap \overline{A}) \cup (\emptyset \cup M) = (N - A) \cup M$.

From here we get the unique solution: $X = (N - A) \cup M$.

We test $f((N-A) \cup M) = (((N-A) \cup M) \cap A, ((N-A) \cup M) \cup A)$

but

$((N-A) \cup M) \cap A = ((N \cap \overline{A}) \cup M) \cap A = ((N \cap \overline{A}) \cap A) \cup (M \cap A) =$
$= ((N \cap (\overline{A} \cap A)) \cup M = (N \cap \emptyset) \cup M = \emptyset \cup M = M$

and

$((N-A) \cup M) \cup A = (N-A) \cup (M \cup A) = (N-A) \cup A = (N \cap \overline{A}) \cup A =$
$= (N \cup A) \cap (\overline{A} \cup A) = N \cap E = N$, $f((N-A) \cup M) = (M, N)$.

Thus $f(P(E)) = K$.

   c. From point a. we have that $g$ is injective, from point b. we have that $g$ surjective, thus $g$ is bijective. The inverse function is:

$g^{-1}(M, N) = (N - A) \cup M$.

**Problem 3.** Let $E \neq \emptyset$, $A, B \in P(E)$ and $f : P(E) \to P(E) \times P(E)$, where $f(X) = (X \cap A, X \cap B)$.
   a. Give the necessary and sufficient condition such that $f$ is injective.
   b. Give the necessary and sufficient condition such that $f$ is surjective.
   c. Supposing that $f$ is bijective, compute its inverse.
*Solution.*

   a. Suppose that $f$ is injective. Then:

$f(A \cup B) = ((A \cup B) \cap A, (A \cup B) \cap B) = (A, B) = (E \cap A, E \cap B) = f(E)$,

from where $A \cup B = E$.

Now we suppose that $A \cup B = E$, it results that:
$X = X \cap E = X \cap (A \cup B) = (X \cap A) \cup (X \cap B) = (Y \cap A) \cup (Y \cap B) = Y \cap (A \cup B) = Y \cap E = Y$

namely from $f(X) = f(Y)$ we obtain that $X = Y$, namely $f$ is injective.

   b. Suppose that $f$ is surjective, for any $M, N \in P(A) \times P(B)$, there exists

$X \in P(E), f(X) = (M, N), (X \cap A, X \cap B) = (M, N), X \cap A = M, X \cap B = N$.

In special cases $(M, N) = (A, \emptyset)$, there exists $X \in P(E)$, from
$X \supset A, \emptyset = X \cap B \supset A \cap B, A \cap B = \emptyset$.

Now we suppose that $A \cap B = \emptyset$ and show that it is surjective.

Let $(M, N) \in P(A) \times P(B)$, then $M \subset A, N \subset B, M \cap B \subset A \cap B = \emptyset$, and $N \cap A \subset B \cap A = \emptyset$, namely $M \cap B = \emptyset$, $N \cap A = \emptyset$ and
$f(M \cup N) = ((M \cup N) \cap A, (M \cup N) \cap B) =$
$= ((M \cap A) \cup (N \cap A), (M \cap B) \cup (N \cap B)) = (M \cup \emptyset, \emptyset \cup N) = (M, N)$,

for any $(M, N)$ there exists $X = M \cup N$ such that $f(X) = (M, N)$, namely $f$ is surjective.

    c. We'll show that $f^{-1}((M, N)) = M \cup N$.

**Remark.** In the previous two problems we can use the characteristic function of the set as in the first problem. We leave this method for the readers.

**Application.** Let $E \neq \emptyset$, $A_k \in P(E)$ $(k = 1,...,n)$ and $f : P(E) \to P^n(E)$, where $f(X) = (X \cap A_1, X \cap A_2, ..., X \cap A_n)$.

Prove that $f$ is injective if and only if $\bigcup_{k=1}^{n} A_k = E$.

**Application.** Let $E \neq \emptyset$, $A_k \in P(E)$, $(k = 1,...,n)$ and $f : P(E) \to P^n(E)$, where $f(X) = (X \cap A_1, X \cap A_2, ..., X \cap A_n)$.

Prove that $f$ is surjective if and only if $\bigcap_{k=1}^{n} \overline{A_k} = \emptyset$.

**Problem 4.** We name the set $M$ convex if for any $x, y \in M$ $tx + (1-t)y \in M$, for any $t \in [0,1]$.

Prove that if $A_k$, $(k = 1,...,n)$ are convex sets, then $\bigcap_{k=1}^{n} A_k$ is also convex.

**Problem 5.** If $A_k$, $(k = 1,...,n)$ are convex sets, then $\bigcap_{k=1}^{n} A_k$ is also convex.

**Problem 6.** Give the necessary and sufficient condition such that if $A$, $B$ are convex/concave sets, then $A \cup B$ is also convex/concave. Generalization for the $\mathbb{N}$ set.

**Problem 7.** Give the necessary and sufficient condition such that if $A$, $B$ are convex/concave sets then $A \Delta B$ is also convex/concave. Generalization for the $\mathbb{N}$ set.

**Problem 8.** Let $f, g : P(E) \to P(E)$, where $f(x) = A - X$, and
$$g(x) = A \Delta X, \ A \in P(E).$$
Prove that $f$, $g$ are bijective and compute their inverse functions.

**Problem 9.** Let $A \circ B = \{(x, y) \in \square \times \square \ | \ \exists z \in \square : (x, z) \in A \text{ and } (z, y) \in B\}$. In a particular case let $A = \{(x, \{x\}) | x \in \square \}$ and $B = \{(\{y\}, y) | y \in \square \}$.
Represent the $A \circ A$, $B \circ A$, $B \circ B$ cases.

**Problem 10.**
    i.    If $A \cup B \cup C = D$, $A \cup B \cup D = C$, $A \cup C \cup D = B$, $B \cup C \cup D = A$, then $A = B = C = D$

ii. Are there different $A$, $B$, $C$, $D$ sets such that
$$A \cup B \cup C = A \cup B \cup D = A \cup C \cup D = B \cup C \cup D?$$

**Problem 11.** Prove that $A \Delta B = A \cup B$ if and only if $A \cap B = \emptyset$.

**Problem 12.** Prove the following identity.
$$\bigcap_{i,j=1, i<j}^{n} A_k \cup A_j = \bigcup_{i=1}^{n} \left( \bigcap_{j=1, j \neq i}^{n} A_j \right)$$

**Problem 13.** Prove the following identities.
$$(A \cup B) - (B \cap C) = (A - (B \cap C)) \cup (B - C) = (A - B) \cup (A - C) \cup (B - C)$$

and

$$A - [(A \cap C) - (A \cap B)] = (A - \bar{B}) \cup (A - C).$$

**Problem 14.** Prove that $A \cup (B \cap C) = (A \cup B) \cap C = (A \cup C) \cap B$ if and only if $A \subset B$ and $A \subset C$.

**Problem 15.** Prove the following identities:
$$(A - B) - C = (A - B) - (C - B),$$
$$(A \cup B) - (A \cup C) = B - (A \cap C),$$
$$(A \cap B) - (A \cap C) = (A \cap B) - C.$$

**Problem 16.** Solve the following system of equations:
$$\begin{cases} A \cup X \cup Y = (A \cup X) \cap (A \cup Y) \\ A \cap X \cap Y = (A \cap X) \cup (A \cap Y) \end{cases}.$$

**Problem 17.** Solve the following system of equations:
$$\begin{cases} A \Delta X \Delta B = A \\ A \Delta Y \Delta B = B \end{cases}.$$

**Problem 18.** Let $X$, $Y$, $Z \subseteq A$. Prove that:
$Z = (X \cap \bar{Z}) \cup (Y \cap \bar{Z}) \cup (\bar{X} \cap Z \cap \bar{Y})$ if and only if $X = Y = \emptyset$.

**Problem 19.** Prove the following identity:
$$\bigcup_{k=1}^{n} [A_k \cup (B_k - C)] = \left( \bigcup_{k=1}^{n} A_k \right) \cup \left[ \left( \bigcup_{k=1}^{n} A_k \right) - C \right].$$

**Problem 20.** Prove that: $A \cup B = (A - B) \cup (B - A) \cup (A \cap B)$.

**Problem 21.** Prove that:

$$(A \Delta B) \Delta C = (A \cap \bar{B} \cap \bar{C}) \cup (\bar{A} \cap B \cap \bar{C}) \cup (\bar{A} \cap \bar{B} \cap C) \cup (A \cap B \cap C).$$